\documentclass[12pt]{amsart}\usepackage{amssymb}

\newtheorem{thm}{Theorem}

\theoremstyle{definition}

\newcommand\Z{\mathbb{Z}}

\newcommand{\heading}[1]{\medskip\noindent\textbf{#1.}}

\begin{document}

\title{Curvature and computation} 
\author{Jon McCammond}
\address{Dept. of Math., University of California, Santa Barbara, CA 93106} 
\thanks{Partial support by the National Science Foundation is
 gratefully acknowledged}
\email{jon.mccammond@math.ucsb.edu}
\date{\today}

\begin{abstract}
  When undergraduates ask me what geometric group theorists study, I
  describe a theorem due to Gromov which relates the groups with an
  intrinsic geometry like that of the hyperbolic plane to those in
  which certain computations can be efficiently carried out.  In
  short, I describe the close but surprising connection between
  negative curvature and efficient computation.  This theorem was one
  of the clearest early indications that applying a metric perspective
  to traditional group theory problems might lead to new and important
  insights.
\end{abstract}

\maketitle

The theorem I want to discuss asserts that there is a close
relationship between two collections of groups: one collection is
defined geometrically and the other is defined computationally.  The
first section describes the relevant geometric and topological ideas,
the second discusses the key algebraic and computational concepts, and
the short final section describes the relationship between them.  An
informal style, similar to the one I use when answering this question
face-to-face, is maintained throughout.

\section{Geometry and topology}

The first thing to highlight is that there is a close relationship
between groups and topological spaces.  More specifically, to each
connected topological space $X$ there is an associated group $G$
called its \emph{fundamental group} and absolutely every group arises
in this way (in the sense that for each group $G$ one can construct a
topological space $X$ whose fundamental group is isomorphic to $G$).
Because of this connection and because spaces with isomorphic
fundamental groups share many key properties, we can use the topology
of the space $X$ to understand the algebraic structure of its
fundamental group $G$.

\heading{\bf Fundamental groups} In order to make this discussion as
accessible as possible, here is a quick sketch of the basic idea
behind the notion of a fundamental group.  As an initial attempt, one
could try to form a group out of a space by using the paths in the
space as our elements and the operation of concatenation as our
multiplication, but there are problems that arise.  First, we want to
be able to ``multiply'' (i.e. concatenate) any two paths, but to do so
we need the first path to end where the second path begins.  To fix
this we select a point $x$ in our space and consider only those paths
that start and end at $x$.  The role of an identity element is played
by the trivial path that starts at $x$ and stays at $x$.

But now we come to the second problem.  Concatenating paths only makes
them longer so that nontrivial paths do not yet have inverses.  It
seems intuitive that traveling along a path in the opposite direction
should count as its inverse but in order to make this work, we need to
replace individual paths with equivalence classes of paths.  We call
two paths equivalent when we can continuously deform one to the other
without moving its endpoints and the multiplication of equivalence
classes of paths is defined as the equivalence class of the
concatenation of representative paths.  It is relatively easy to check
that this multiplication is well-defined and that the resulting
algebraic structure is a group.  Moreover, so long as the space $X$ is
path connected, the algebraic structure of this group does not depend
on our choice of basepoint up to isomorphism.  In other words, the
group $G$ is an invariant of the space $X$ itself independent of our
choice of basepoint.

The standard illustration of this procedure and in many ways the most
crucial one is the following: the fundamental group of the unit circle
is isomorphic to the integers.  The nontrivial group elements come
from paths that wrap around the circle and, in fact, the equivalence
classes essentially collect together those paths that wrap around the
circle the same number of times in the same direction.  A second
example, closely related to the first, is that the fundamental group
of the torus is $\Z \oplus \Z$.

\heading{\bf Covering spaces} There is another connection between the
space $X$ and its fundamental group $G$.  So long as the space $X$ is
sufficiently nice, it can be completely unwrapped in the following
sense.  There is another space called its \emph{universal cover} with
trivial fundamental group and a projection map back to $X$ that is
locally a homeomorphism.  In the case of the circle, its universal
cover is an infinite spiral, continuing forever in both directions.
Topologically this space looks like the real line, it is contractible
and it is easy to believe that its fundamental group is trivial.  It
is also easy to see that the natural projection from the spiral to the
circle is a local homeomorphism.  The universal cover of the torus is
the euclidean plane and the projection map from the plane to the torus
is the one that first wraps it up into an infinite cylinder in one
direction and then wraps it up in the other direction into a torus.

One useful fact is that the fundamental group of $X$ acts on its
universal cover by homeomorphisms.  In our examples, the integers act
on the infinite spiral by rigidly shifting it up or down and the group
$\Z \oplus \Z$ acts on the plane by translating by vectors with
integer coordinates.  In fact, the action of the fundamental group on
the universal cover is always transitive on the preimages of a point
and these preimages are in one-to-one correspondence with the elements
of the fundamental group.  From the early twentieth century to the
present day group theorists have used this relationship to study
infinite discrete groups.  In particular, the topology of the
universal cover on which the fundamental group $G$ acts can be used to
extract information about the algebraic structure of $G$.

\heading{\bf Gestures} This is probably as good as place as any to
mention one key aspect of my interactional style that is difficult to
replicate in a written text.  Gestures are an important aspect of how
I communicate mathematics orally and this is especially true when I am
talking to students.  In particular, throughout this entire discussion
I usually employ a collection of gestures in specific locations to
focus attention and to illustrate what is going on.  The result is
something like Prokofiev's orchestration that accompanies the story of
Peter and the Wolf.  Recurring characters (mathematical concepts) have
musical themes (stylized gestures) that are repeated every time they
reappear.  Whenever I mention the geometric and topological aspects of
groups I gesture to my lefthand side (and the algebraic and
computational aspects involve gestures to my righthand side).  The
space $X$ is located on the lower left and its universal cover is
directly above it.  The gesture associated to $X$, its theme, is the
miming of the shape of a torus.  For its universal cover, I start with
the torus on the lower left and then raise my arms and spread out my
hands to indicate the euclidean plane.  The action of $G$ on the
universal cover of $X$ is indicated by moving both hands (in euclidean
plane position) in small syncronized circles with the hands themselves
always pointing in the same direction and maintaining a rigid
relationship between them.  The reader might want to visualize these
gestures as they read along.

\heading{\bf Metrics on groups} Returning to our discussion of the
action of the fundamental group on the universal cover, suppose we add
a metric to the original space $X$.  In our example, rather than
imagining a space that is merely a topological torus, imagine a space
with a precise metric so that we can calculate distances, angles and
areas.  This local metric information induces a metric on the
unwrapped version and we can use this metric on the universal cover to
turn the group $G$ itself into a metric space.  I should also point
out that the action of $G$ on the universal cover is one that
preserves this local metric information.  In other words, it acts on
the universal cover by isometries.

To turn $G$ into a metric space we to pick a point $\tilde x$ in the
universal cover and then for each group element $g \in G$ record the
distance between this point and its image under $g$.  We call this the
distance from the identity element to $g$.  More generally, given two
group elements $g$ and $g'$, we define the distance between them to be
the distance in the universal cover between the images of $\tilde x$
under $g$ and $g'$.  It is now relatively easy to convince yourself
that this distance function defines a metric on the elements of $G$,
i.e. it is symmetric, nonzero on distinct pairs of elements and
satisfies the triangle inequality.

\heading{\bf Intrinsic Metrics} The other thing that is obvious is
that the precise values of the metric on $G$ very much depend on the
specific metric we added to our space $X$ and the point $\tilde x$
that we selected.  It turns out, however, that when $X$ satisfies
certain minimal conditions (such as being compact) altering the metric
on $X$ or choosing a different point $\tilde x$ does not significantly
alter the induced metric on $G$.  More specifically, given two
distinct metrics on $X$ and two different selected points, the metrics
they induce on $G$ are related by linear inequalities.  In other
words, there exist constants so that for every pair of elements in
$G$, their distance in the first metric is bounded above by a linear
function of their distance in the second metric and vice versa.  Two
metrics that are related in this way are said to be quasi-isometric
and the notion of quasi-isometry partitions all metrics into
quasi-isometry classes.  In this language, the result I'm alluding to
is that for any reasonable space $X$ with fundamental group $G$, the
possible metrics on $X$ induce metrics on $G$ that all belong to the
same quasi-isometry class.  In fact, this remains true even if we
replace our reasonable space $X$ with any other reasonable space $Y$
with the same fundamental group $G$, a result known as the
Milnor-Svarc Theorem.  This means that if $G$ is the kind of group
that can be the fundamental group of a reasonable metric space $X$,
then the quasi-isometry class of the metric induced on $G$ through its
action on the universal cover of $X$ is completely independent of the
space $X$ used to produce this metric.  This is what geometric group
theorists mean when they say that (reasonable) groups come equipped
with an intrinsic metric that is well-defined up to quasi-isometry.

\heading{\bf Differential geometry} We are now going to restrict our
attention to a special class of groups but the motivation for this
restriction involves a short digression into the history of a
different part of mathematics.  Differential geometry is an area that
studies spaces called Riemannian manifolds that are locally
homeomorphic to $n$-dimensional euclidean space and which come
equipped with a nice smooth metric that allows them to be investigated
using the standard tools of multivariable calculus.  Early on,
differential geometers defined various notions of curvature and they
proved that Riemannian manifolds that are negatively curved, in a
suitable sense, have very nice properties such as a contractible
universal cover.  Initially their proofs used the full force of the
analytic tools available to them, but as they simplified the proofs to
extract the essence of why these results were true, they soon
discovered that they could assume much less about the original space
and still produce significant consequences.  In fact, all that was
really necessary was that certain inequalities hold involving points
on the sides of geodesic triangles.  Once reformulated in this way,
their ideas could be applied to a much larger class of metric spaces
which did not necessarily locally look like euclidean space and where
the ordinary operations of multivariable calculus could not be
applied.  One of these differential geometers was Misha Gromov and he
soon realized that these distilled ideas from differential geometry
could be applied to infinite discrete groups.

\heading{\bf Thin triangles} The key definition is inspired by the
properties of triangles in the hyperbolic plane.  If you have ever
studied the geometry of the hyperbolic plane, you have probably
learned that there are important differences between triangles in the
hyperbolic plane and triangles in the euclidean plane.  In a euclidean
triangle, the sum of its three angles is $\pi$ but in a hyperbolic
triangle, the sum of its angles is always strictly less than $\pi$.
The more relevant fact about hyperbolic triangles for our discussion
is one that does not always make it into a first course on hyperbolic
geometry, namely, that all triangles in the hyperbolic plane are
uniformly thin.

In the euclidean plane, some triangles are fat.  What I mean by this
is that for every constant $r$ we can find a euclidean triangle and a
point in its interior so that the distance from this point to any
point on its boundary is at least $r$.  For example, the center of a
large equilateral triangle has this property.  In the hyperbolic plane
you can do this for small values of $r$ but not for large values of
$r$.  Let me give an equivalent reformulation of this property where I
actually know the exact value of the constant where the behavior
changes.  In the euclidean plane for any constant $r$ it is easy to
find a triangle and a point $p$ on one of its sides so that the
distance from $p$ to any point on either of the other two sides is at
least $r$.  In the hyperbolic plane it turns out that given any
triangle and any point $p$ on one of its sides, there is a point $q$
on one of its other sides so that the distance from $p$ to $q$ is less
than $\log(1 + \sqrt{2})$.  This exact value is relatively easy to
establish but the interesting point is that such a value even exists.

\heading{\bf Hyperbolicity} Gromov turned the uniform thinness of
triangles in the hyperbolic plane into a defining characteristic of
hyperbolic spaces and hyperbolic groups.  A space is called
\emph{$\delta$-hyperbolic} when all geodesic triangles in this space
are $\delta$-thin for a fixed constant $\delta$.  In other words,
given any three points and any three length-minimizing paths
connecting them into a triangle and given any point $p$ on one of
these paths, there is a point $q$ on one of the other two paths so
that the distance from $p$ to $q$ is less than $\delta$.  A group $G$
is called \emph{word hyperbolic} or \emph{Gromov hyperbolic} when it
is the fundamental group of a reasonable metric space $X$ whose
universal cover is $\delta$-hyperbolic for some constant $\delta$.  It
turns out that being Gromov hyperbolic really is an intrinsic property
of $G$ in the sense that it is independent of our choice of $X$ and of
our choice of a metric on $X$.  Concretely, if $X$ and $Y$ are
reasonable metric spaces with fundamental group $G$ and one of them
has a $\delta$-hyperbolic universal cover then the other universal
cover is $\delta'$-hyperbolic for a possibly different constant
$\delta'$.  In short, groups that are hyperbolic in the sense of
Gromov are those where the geometry of its intrinsic metric shares a
key property possessed by triangles in the hyperbolic plane.

\section{Algebra and computation}

And now for something completely different.  Set the geometric and
topological properties of groups aside for the moment and consider
their algebraic and computational properties.  The first thing to note
is that the infinite groups which are the easiest to work with from a
computational perspective are those that have some sort of finite
description.

\heading{\bf Descriptions of groups} The classical method of
describing an infinite group is to list a set of elements that are
sufficient to generate the entire group and then to list some
relations satisfied by these elements that are sufficient to generate
all of the relations that hold in the group.  Such a
\emph{presentation} is said be finite when both the set of generators
and the set of relations are finite and the group it describes is
called \emph{finitely presented}.  The classical example of a finite
presentation is the group generated by $a$ and $b$ and subject only to
the relation that $ab = ba$.  This is a finite description of the
group $\Z \oplus \Z$.  There are other ways to characterize the class
of finitely presented groups that make clear that this is an important
and interesting class of groups to study.  For example, the class of
finitely presented groups is exactly the same as the class of groups
that are fundamental groups of compact manifolds and it is exactly the
same as the class of groups that are fundamental groups of finite
simplicial complexes.

\heading{\bf The word problem} It is traditional to use a language
metaphor when working with a finitely presented group.  Individual
generators are called \emph{letters} and finite products of generators
and their inverses are called \emph{words}.  One problem that
immediately arises is that because the generators satisfy relations,
there are typically many different words that represent the exact same
element of the group.  In the standard presentation of the group $\Z
\oplus \Z$, for example, both $aaabb$ and $ababa$ represent the same
element even though they are distinct words.  The key question, first
identified by Max Dehn in 1912, is the \emph{word problem}: For a
fixed finite presentation, is there an algorithm that takes as input
two words written as products of the generators and their inverses and
outputs whether or not they represent the same element of the group
after a finite amount of time.  For $\Z \oplus \Z$ the answer is yes,
there does exist such an algorithm.  One such algorithm goes as
follows.  Systematically move all the $a$'s and $a^{-1}$'s to the left
and all the $b$'s and $b^{-1}$'s to the right and then simplify until
the final result is a word of the form $a^i b^j$ for some integers $i$
and $j$.  Two words that have the same \emph{normal form} represent
the same group element and two words that have distinct normal forms
represent distinct group elements.  This works for $\Z \oplus \Z$ but
the general situation is much more complicated.

\heading{\bf Some problems cannot be solved} In the early twentieth
century mathematicians were beginning to learn that there is an
important distinction between what is true and what can be proved.  In
the same way that a statement such as ``This is a lie'' cannot
consistently be assigned a truth value, G\"odel showed how one could
construct a problematic assertion in any finite axiomatic system for
the natural numbers.  This problematic assertion is either a true
statement that cannot be proved from the axioms, or it is a false
statement that the axioms can prove.  This means that in any
consistent axiom system for the natural numbers there are things that
are true but not provable.  When translated into the language of the
theory of computation, this means that there are problems that cannot
be solved algorithmically and one can prove that they cannot be solved
algorithmically.  One example of such an unsolvable problem is the
halting problem: Does there exist a computer program which takes as
input an arbitrary computer program and outputs, after a finite amount
of time, whether or not the program given as input will run forever?
The answer is that no such generic program analyzing software can
exist since there will always be some program that it cannot
successfully analyze.

Once mathematicians realized that some problems cannot be solved, they
used that fact to prove that other problems cannot be solved.  They
did this by showing that a solution to the second problem leads to a
solution to the first problem, which contradicts the fact that we know
the first problem cannot be solved.  Working along these lines, Boone
and Novikov, working independently, showed that a single algorthm that
solves the word problem in an arbitrary finitely presented group
cannot and does not exist.  In fact, there are explicit finite
presentations for which it is known that there is no algorithm to
solve the word problem for this specific group.

\heading{\bf Efficient solutions} The fact that some finitely
presented groups have word problems that cannot be solved merely
prompts mathematicians to shift their attention to those groups with
word problems that can be solved.  Going one step further, we can
divide groups with solvable word problems into classes based on how
hard their word problems are to solve.  One indication of the level of
difficulty is how long it takes for the algorithm to work: how many
units of time does it take as a function of the total length of the
two words under consideration?  In other words, let $n$ denote this
total length of the two words given as input and describe the time
bound as a function of $n$.  Is it a linear function of $n$?
quadratic?  polynomial?  exponential?

Clearly the best possible algorithms cannot run in less than linear
time since in order to correctly answer whether or not two words
represent the same element in the group, the program must, at the very
least, read the two words which takes a linear amount of time.  In our
example of the standard presentation of $\Z \oplus \Z$, the process
described that places words in normal form $a^i b^j$ can take a
quadratic amount of time since one of the initial inputs might be $b^m
a^m$.  Each $b$ must be moved past each $a$ which involves $m^2$ local
modifications.  These can be visualized in the plane as pushing across
single squares from one pair of sides of a large rectangle to the
other pair of sides.

\heading{\bf Computational complexity} At this point I need to make a
slight technical aside in order to make my later descriptions precise.
There are various standard models of computation that one can use.
Some notions of computational complexity (such as polynomial time) are
quite robust in the sense that they define the same class of problems
regardless of the model one uses, but linear time is not one of them.
There are, for example, problems which one can solve in linear time on
a multi-tape Turing machine that take longer on a single-tape Turing
machine.  Concretely, one can solve the word problem for $\Z \oplus
\Z$ in linear time by reading the words $u$ and $v$ and merely
checking that they have the same number of $a$'s and $b$'s after
cancelation, but this is not the type of algorithm that I am
interested in.  When discussing groups that have a linear time
solution to their word problem, I have in mind a very specific type of
algorithm.  The algorithms I wish to consider are those that work by
taking the words and systematically rewriting them using the relations
in the presentation.  One might call these \emph{relation-driven
  algorithms} for the word problem.

\heading{\bf Dehn's algorithm} There is a famous relation-driven
algorithm called \emph{Dehn's algorithm} which is easy to implement
even though it only works for special presentations of certain groups.
It proceeds as follows.  Start reading a word such as $uv^{-1}$ from
the beginning and look for a subword that represents strictly more
than half of one of the relations.  If one is found, replace it with
the shorter half of the relation and back up to the beginning of the
replacement.  Continue.  For the presentations where Dehn's algorithm
works, this procedure terminates in linear time and it produces the
trivial word if and only if the input word is equal to the identity in
the group.  This relation-driven algorithm fails when there is a word
equivalent to the identity of the group which does not contain more
than half of a relation.  The standard presentation of the fundamental
group of an orientable surface of genus at least $2$ (i.e. with more
than one hole) is one where Dehn's algorithm is known to work.

\heading{\bf Isoperimetric inequalities} The best possible time bound
for a relation-driven algorithm solution to the word problem in a
particular finitely presented group $G$ is closely related to the
isoperimetric inequality satisfied by closed loops in the universal
cover of any reasonable metric space $X$ with fundamental group $G$.
By isoperimetric inequality we mean the following.  For each closed
loop in a simply connected metric space, we can measure the minimal
area of a disc mapped isometrically into the space so that its
boundary is the specified closed loop.  We measure the ratio of this
minimal area to the length of the curve and then find the largest such
ratio as the curves vary over all curves up to a specified length.  As
this bound grows we get an increasing function that measures how hard
it is to fill loops of a given bounded size.  It turns out that the
rate of growth of this function is independent of $X$ and is an
invariant of $G$ alone.  It is bounded above by a recursive function
iff the word problem for $G$ can be solved and for our example of $\Z
\oplus \Z$ the growth rate of this function is quadratic.  For
finitely presented groups, the isoperimetric inequality essentially
measures the number of relations needed to prove a particular word is
equal to the identity, and thus it provides a lower bound on the time
it takes for a relation-driven algorithm to solve the word problem.
In particular, the finitely presented groups whose word problem can be
solved in linear time by a relation-driven algorithm must have an
isoperimetric inequality that grows linearly.  These are groups with
the fastest possible solution to their word problem, the groups that
arguably have the best possible computational properties.

\section{Gromov's theorem}

And now for the surprising connection.  In the 1980s Gromov proved
that the groups with the best possible computational properties
correspond \emph{exactly} to those with an intrinsic geometry that is
negatively curved in the sense described earlier \cite{Gr85}.

\begin{thm}
  A finitely presented group has a linear time solution to its word
  problem (in the sense described above) \emph{if and only if} it is
  hyperbolic in the sense of Gromov. 
\end{thm}

A more technical version of this theorem would be that for a fixed
finitely presented group the following are equivalent: (1) it has an
alternative presentation where Dehn's algorithm works, (2) there is a
relation-driven algorithm that solves the word problem in linear time,
(3) it has a linear isoperimetric inequality, and (4) it is hyperbolic
is the sense of Gromov.  We have already mentioned that (1) implies
(2) implies (3).  Gromov proved that (3) implies (4) implies (1).

This is an amazing result because the implications are not merely in
one direction; it is an exact correspondence.  As a consequence of
this result, geometric group theorists tend to view hyperbolic
geometry as the best possible geometry for a group to have since it
corresponds to the group having the best possible computational
properties.  This is the kind of result that makes researchers sit up
and take notice, and it prompted a thorough-going review of the
foundations of the subject.  It also immediately prompts a large
number of follow-up questions.  How strong is this bridge between
geometry and topology on the one hand and algebra and computation on
the other? In particular, what happens when we expand the class of
groups under consideration?  Are there geometric consequences when a
finitely presented group has a quadratic, cubic or polynomial time
solution to its word problem?  Are there computational consequences
when a finitely presented group has an intrinsic geometry that is
non-positively curved in some sense (rather than negatively curved)?
Over the past 30 years these types of questions have led to the
development of several general theories such as the theory of
automatic and biautomatic groups \cite{ECHLPT92} and the theory of
groups that act on nonpositively-curved spaces \cite{BrHa99}.

A second natural set of questions involves taking the various classes
of groups traditionally investigated by combinatorial group theorists
(such as outer automorphisms of free groups, mapping class groups of
closed surfaces, braid groups, Coxeter groups, Artin groups,
one-relator groups, etc., etc.) and asking which of the various
general theories of curvature and computation apply in each case.  To
my mind, this single theorem of Gromov is like the Big Bang and it
played a major role in the creation of a new subfield called geometric
group theory.

\newcommand{\etalchar}[1]{$^{#1}$}
\providecommand{\bysame}{\leavevmode\hbox to3em{\hrulefill}\thinspace}
\providecommand{\MR}{\relax\ifhmode\unskip\space\fi MR }
\providecommand{\MRhref}[2]{%
  \href{http://www.ams.org/mathscinet-getitem?mr=#1}{#2}
}
\providecommand{\href}[2]{#2}

\end{document}